\documentclass[a4paper, 11pt]{amsart}
\usepackage{graphicx}
\usepackage{bm}
\usepackage{amssymb}

\begin{document}

\title{Group geometrical axioms for magic states of quantum computing}

\author{Michel Planat$\dag$, Raymond Aschheim$\ddag$,\\ Marcelo M. Amaral$\ddag$ and Klee Irwin$\ddag$}

\address{$\dag$ Universit\'e de Bourgogne/Franche-Comt\'e, Institut FEMTO-ST CNRS UMR 6174, 15 B Avenue des Montboucons, F-25044 Besan\c con, France.}
\email{michel.planat@femto-st.fr}

\address{$\ddag$ Quantum Gravity Research, Los Angeles, CA 90290, USA}
\email{raymond@QuantumGravityResearch.org}
\email{Klee@quantumgravityresearch.org}
\email{Marcelo@quantumgravityresearch.org}

\begin{abstract}

Let $H$ be a non trivial subgroup of index $d$ of a free group $G$ and $N$ the normal closure of $H$ in $G$. The coset organization in a subgroup $H$ of $G$ provides a group $P$ of permutation gates whose common eigenstates are either stabilizer states of the Pauli group or magic states for universal quantum computing. A subset of magic states consists of MIC states associated to minimal informationally complete measurements. It is shown that, in most cases, the existence of a MIC state entails that the two conditions (i) $N=G$ and (ii) no geometry (a triple of cosets cannot produce equal pairwise stabilizer subgroups), or that these conditions are both not satisfied.  Our claim is verified by defining the low dimensional MIC states from subgroups of the fundamental group $G=\pi_1(M)$ of some manifolds encountered in our recent papers, e.g. the $3$-manifolds attached to the trefoil knot and the figure-eight knot, and the $4$-manifolds defined by $0$-surgery of them. 
Exceptions to the aforementioned rule are classified in terms of geometric contextuality (which occurs when cosets on a line of the geometry do not all mutually commute).


 \end{abstract}

\maketitle


\normalsize

\section{Introduction}

Interpreting quantum theory is a long standing effort and not a single approach can exhaust all facets of this fascinating subject. Quantum information owes much to the concept of a (generalized) Pauli group for understanding quantum observables, their commutation, entanglement, contextuality and many other aspects, e.g. quantum computing. Quite recently, it has been shown that quantum commutation relies on some finite geometries such as generalized polygons and polar spaces \cite{Pauli2011}. Finite geometries connect to the classification of simple groups as understood by prominent researchers as Jacques Tits, Cohen Thas and many others \cite{Thas2005, Zoology2017}.

In the Atlas of finite group representations \cite{Atlasv3}, one starts with a free group $G$ with relations, then the finite group under investigation $P$ is the permutation representation of the cosets of a subgroup of finite index $d$ of $G$ (obtained thanks to the Todd-Coxeter algorithm). As a way of illustrating this topic, one can refer to \cite[Table 3]{Zoology2017} to observe that a certain subgroup of index $15$ of the symplectic group $S_4'(2)$ corresponds to the $2QB$ (two-qubit) commutation of the $15$ observables in terms of the generalized quadrangle of order two, denoted $GQ(2,2)$ (alias the doily). For $3QB$, a subgroup of index $63$ in the symplectic group $S_6(2)$ does the job and the commutation relies on the symplectic polar space $W_5(2)$ \cite[Table 7]{Zoology2017}. An alternative way to approach $3QB$ commutation is in terms of the generalized hexagon $GH(2,2)$ (or its dual) which occurs from a subgroup of index $63$ in the unitary group $U_3(3)$ \cite[Table 8]{Zoology2017}. Similar geometries can be defined for multiple qudits (instead of qubits).

The straightforward relationship of quantum commutation to the appropriate symmetries and finite groups was made possible thanks to techniques developed by the first author (and coauthors) that we briefly summarize. This will be also useful at a later stage of the paper with the topic of magic state quantum computing. 

The rest of this introduction recalls how the permutation group organizing the cosets leads to the finite geometries of quantum commutation (in Sec. \ref{Int1}) and how it allows the computation of magic states of universal quantum commutation (in Sec. \ref{Int2}). In this paper, it is shown that magic states themselves may be classified according to their coset geometry with two simple axioms (in Sec. \ref{Int3}).

\subsection{Finite geometries from cosets \cite{Zoology2017,PGHS2015,Planat2015}}
\label{Int1}
One needs to define the rank $r$ of a permutation group $P$. First~it is expected that $P$ acts faithfully and transitively on the set $\Omega=\{1,2,\cdots, n\}$ as a subgroup of the symmetric group $S_n$. The action of $P$ on a pair of distinct elements of $\Omega$ is defined as $(\alpha,\beta)^p=(\alpha^p,\beta^p)$, $p\in P$, $\alpha \ne \beta$. The orbits of $P$ on $\Omega \times \Omega$ are called orbitals, and the number of orbits is called the rank $r$ of $P$ on $\Omega$. The rank of $P$ is at least two, and the two-transitive groups identify to the rank two permutation groups.
 Next, selecting a pair $(\alpha,\beta)\in \Omega \times \Omega$, $\alpha \ne \beta$, one introduces the two-point stabilizer subgroup $P_{(\alpha,\beta)}=\{p \in P|(\alpha,\beta)^p=(\alpha,\beta)\}$. There exist $1 < m \le r$ such non-isomorphic (two-point stabilizer) subgroups $S_m$ of $P$.
 Selecting one with $\alpha \ne \beta$, one defines a point/line incidence geometry $\mathcal{G}$ whose points are the elements of $\Omega$ and whose lines are defined by the subsets of $\Omega$ sharing the same two-point stabilizer subgroup. Thus, two lines of $\mathcal{G}$ are distinguished by their (isomorphic) stabilizers acting on distinct subsets of $\Omega$. A non-trivial geometry arises from $P$ as soon as the rank of the representation $\mathcal{P}$ of $P$ is $r>2$, and simultaneously, the number of non isomorphic two-point stabilizers of $\mathcal{P}$ is $m>2$. 
Further, $\mathcal{G}$ is said to be {\it contextual} (shows {\it  geometrical contextuality}) if at least one of its lines/edges corresponds to a set/pair of vertices encoded by non-commuting cosets \cite{Planat2015}.

Figure 1 illustrates the application of the two-point stabilizer subgroup approach just described for the index $15$ subgroup of the symplectic group is $S'_4(2)=A_6$ whose finite representation is \newline
$H~=~\left\langle a,b|a^2=b^4=(ab)^5=(ab^2)^5=1\right\rangle.$
 The finite geometry organizing the coset representatives is the generalized quadrangle $GQ(2,2)$. The other reresentation is in terms of the two-qubit Pauli operators, as first found in \cite{Pauli2011, SanigaPlanat2007}.
It is easy to check that all lines not passing through the coset $e$ contains some mutually not commuting cosets so that the $GQ(2,2)$ geometry is contextual. The embedded $(3 \times 3)$-grid shown in bold (the so-called Mermin square) allows a $2QB$ proof of Kochen-Specker theorem \cite{Planat2012}.

\begin{figure}[ht]
\includegraphics[width=7cm]{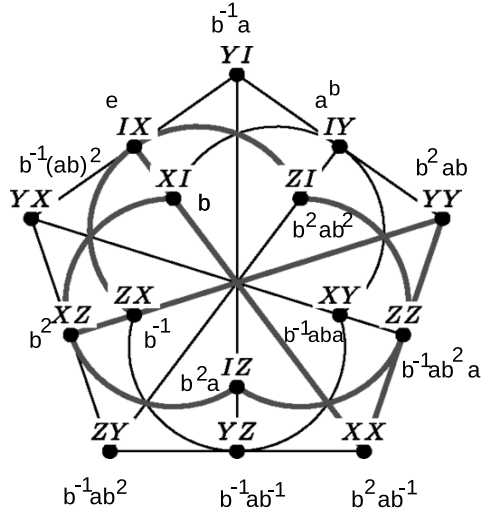}
\caption{The generalized quadrangle of order two $GQ(2,2)$. The picture provides a representation in terms of the fifteen $2QB$ observables that are commuting by triples: the lines of the geometry. Bold lines are for an embedded $3\times 3$ grid (also called Mermin square) that is a basic model of Kochen-Specker theorem (e.g. \cite[Fig.1]{Zoology2017} or \cite{Planat2012}). The second representation is in terms of the cosets of the permutation group arising from the index $15$ subgroup of $G \cong A_6$ (the $6$-letter alternating group).}
\label{fig1}
\end{figure} 

\subsection{Magic states in quantum computing}
\label{Int2}

Now we recall our recent work about the relation of coset theory to the magic states of universal quantum computing.
Bravyi $\&$ Kitaev introduced the principle of \lq magic state distillation' \cite{Bravyi2004}: universal quantum computation, the possibility of getting an arbitrary quantum gate, may be realized thanks to stabilizer operations (Clifford group unitaries, preparations and measurements) and an appropriate single qubit non-stabilizer state, called a \lq magic state'. Then, irrespective of the dimension of the Hilbert space where the quantum states live, a non-stabilizer pure state was called a magic state \cite{Veitch2014}. An improvement of this concept was carried out in \cite{PlanatGedik} showing that a magic state could be at the same time a fiducial state for the construction of a minimal informationally complete positive operator-valued measure, or MIC, under the action on it of the Pauli group of the corresponding dimension. Thus UQC in this view happens to be relevant both to magic states and to MICs. In \cite{PlanatGedik}, a $d$-dimensional magic state is obtained from the permutation group that organizes the cosets of a subgroup $H$ of index $d$ of a two-generator free group $G$. This is due to the fact that a permutation may be realized as a permutation matrix/gate and that mutually commuting matrices share eigenstates - they are either of the stabilizer type (as elements of the Pauli group) or of the magic type. It is enough to keep magic states that are simultaneously fiducial states for a MIC because the other magic states may loose the information carried during the computation. A catalog of the magic states relevant to UQC and MICs can be obtained by selecting $G$ as the two-letter representation of the modular group $\Gamma=PSL(2,\mathbb{Z})$ \cite{PlanatModular}. The next step, developed in \cite{MPGQR1,MPGQR3}, is to relate the choice of the starting group $G$ to three-dimensional topology. More precisely, $G$ is taken as the fundamental group $\pi_1(S^3 \setminus L)$ of a $3$-manifold $M^3$ defined as the complement of a knot or link $L$ in the $3$-sphere $S^3$. A branched covering of degree $d$ over the selected $M^3$ has a fundamental group corresponding to a subgroup of index $d$ of $\pi_1(M^3)$ and may be identified as a sub-manifold of $M^3$, the one leading to a MIC is a model of UQC. In the specific case of $\Gamma$, the knot involved is the left-handed trefoil knot $T_1=3^1$, as shown in \cite{PlanatModular} and \cite[Sec. 2]{MPGQR1}.

\subsection{Coset geometry of magic states}
\label{Int3}

The goal of the paper is to classify the magic states according to the coset geometry where they arise. We start from a non trivial subgroup $H$ of index $d$ of a free group $G$ and  we denote $N$ the normal closure of $H$ in $G$. As above, the coset organization in a subgroup $H$ of $G$ provides a group $P$ of permutation gates whose common eigenstates are either stabilizer states of the Pauli group or magic states for universal quantum computing. A subset of magic states consists of MIC states associated to minimal informationally complete measurements. 

It is shown in the paper  that, in many cases, the existence of a MIC state entails that the two conditions (i) $N=G$ and (ii) no geometry (a triple of cosets cannot produce equal pairwise stabilizer subgroups), or that these conditions are both not satisfied.  Our claim is verified by defining the low dimensional MIC states from subgroups of the fundamental group $G=\pi_1(M)$ of manifolds encountered in our recent papers, e.g. the $3$-manifolds attached to the trefoil knot and the figure-eight knot, and the $4$-manifolds defined by $0$-surgery of them. 

Exceptions to the aforementioned rule are classified in terms of geometric contextuality (which occurs when cosets on a line of the geometry do not all mutually commute).

In section \ref{Fig8}, one deals with the case of MIC states obtained from the subgroups of the fundamental group of Figure-of-Eight knot hyperbolic manifold and its $0$-surgery. In section \ref{TrefoilKnot}, the MIC states produced with the trefoil knot manifold and its $0$-surgery are investigated.

\section{MIC states pertaining to the Figure-of-Eight knot and its $0$-surgery}
\label{Fig8}

We first investigate the relation of MIC states to the group geometrical axioms (i)-(ii) (or their negation) in the context of the Figure-of-Eight knot $K4a1$ (in Sec. \ref{Fig8b}) and its $0$-surgery (in Sec. \ref{Fig8a}). 
The fundamental group of the complement of $K4a1$ in the $3$-sphere $G=\pi_1(S^3 \setminus K4a1)$, and its connection to MICs, is first studied in \cite[Table 2]{MPGQR1}, below are new results and corrections.

\subsection{Group geometrical axioms applied to the fundamental group $\bm{\pi_1(Y)}$}
\label{Fig8a}

The manifold $Y$ defined by $0$-surgery on the knot $K4a1$ is of special interest as shown in \cite[Sec 2]{MPGQR3} and references therein. The number of subgroups of index $d$ of the fundamental group $\pi_1(Y)$ is as follows

\small
$\eta_d[\pi_1(Y)]=[1,1,1, {\bf 2},2,~{\bf 5},1,2,{\bf 2},4,~{\bf 3},17,1,1,2,~3, 1, 6, 3, 6,~1, 3, 1, 43,~\ldots ],$
\normalsize

\noindent
where a bold number means that a MIC exists at the corresponding index.

In Table \ref{tab1}, one summarizes the check of our axioms (i) and (ii) applied to $\pi_1(Y)$. A triangle $\Delta$ means that a geometry does exist (corresponding to at least a triple of cosets with equal pairwise stabilizer subgroups), thus with (ii) is violated. According to our theory, for a MIC to exist, we should have (i) and (ii) satisfied, or both of them violated. The former case occurs for $d=9$, $11$ and $19$. The latter case occurs for $d=6$ where the geometry is that of the octahedron (with the $3$-partite graph $K(2,2,2)$) and $d=20$, where the geometry is encoded by the the complement of the line graph of the bipartite $K(4,5)$. In all of these five cases, a $pp$-valued MIC does exist. 

\begin{figure}[h]
\includegraphics[width=4cm]{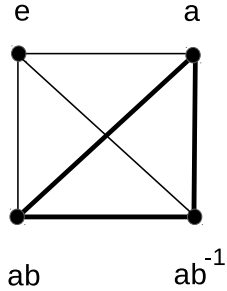}
\caption{The contextual geometry associated to the 2QB-MIC and permutation group $P=A_4$ in Table \ref{tab1}. The line/triangle $\{a,ab,ab^{-1}\}$ is not made of mutually commuting cosets, thus geometric contextuality occurs.}
\label{fig2}
\end{figure}

For dimension $4$, the bold triangle points out a violation since (i) is true and (ii) is false while the $2QB$-MIC exists. In this case the geometry is the tetrahedron (with complete graph $K_4$) but not all cosets on a line/triangle are mutually commuting, a symptom of geometric contextuality, as shown in Fig. \ref{fig2}.


%
\begin{table}[h]
\begin{center}
\begin{tabular}{|c|l|c|c||r|l|r|}
\hline 
\hline
d & P & (i) & pp & geometry\\
\hline
4  & $A_4$&yes &2 & 2QB MIC, $\bm{\Delta}$\\
5  & $10$ &yes & & $\Delta$\\
6  & $A_4$ &no & 2& $6$-dit MIC, $K(2,2,2)$\\
9  &$(36,9) \cong 3^2 \rtimes 4$ &yes & 2& 2QT MIC\\
11  &$(55,1)= 11 \rtimes 5$, $(\times 2)$ &yes & 3& $11$-dit MIC\\
16  & $(48,3) \cong 4 \rtimes A_4 $ &yes & & $\Delta$\\
19  & $(171,3) \cong 19 \rtimes 9$ &yes & 3& $19$-dit MIC\\
20  & $(120,39) \cong 4 \rtimes(5 \rtimes (6,2))$ &yes & & $\overline{\mathcal{L}(K(4,5))}$\\
\hline
\hline
\end{tabular}
\caption{Table of subgroups of the fundamental group $\pi_1[S^3\setminus K4a1(0,1)]$ [with $K4a1(0,1)$ the 0-surgery over the Figure-of-Eight knot].  The permutation group $P$ organizing the cosets in column 2. If (i) is true, unless otherwise specified, the graph of cosets leading to a MIC is that of the $d$-simplex [ and/or the condition (ii) is true: no geometry]. The symbol $\Delta$ means that (ii) fails to be satisfied. When there exists a MIC with (i) true and (ii) false, the geometry is shown in bold characters (here this occurs in dimension $4$, see Fig. \ref{fig2}). If it exists, the MIC is $pp$-valued as given in column 4. In addition, $K(2,2,2)$ is the binary tripartite graph (alias the octahedron) and $\overline{\mathcal{L}(K(4,5))}$ means the complement of the line graph of the bipartite graph $K(4,5)$.} 
\label{tab1}
\end{center}
\end{table}
%

 
\subsection{Group geometrical axioms applied to the fundamental group $\bm{\pi_1(S^3 \setminus K4a1)}$}
\label{Fig8b}

The submanifolds obtained from the subgroups of index $d$ of the fundamental group $\pi_1(S^3 \setminus K4a1)$ for the Figure-of-Eight knot complement are given in Table \ref{tab2} (column 3), as identified in SnapPy \cite{SnapPy} (this corrects a few mistakes of \cite[Table 2]{MPGQR1}). 

As for the subsection above, when axioms (i) and (ii) are simultaneously satisfied (or both are not satisfied), a MIC is created. Otherwise no MIC exist in the corresponding dimension, as expected.

\begin{figure}[h]
\includegraphics[width=5cm]{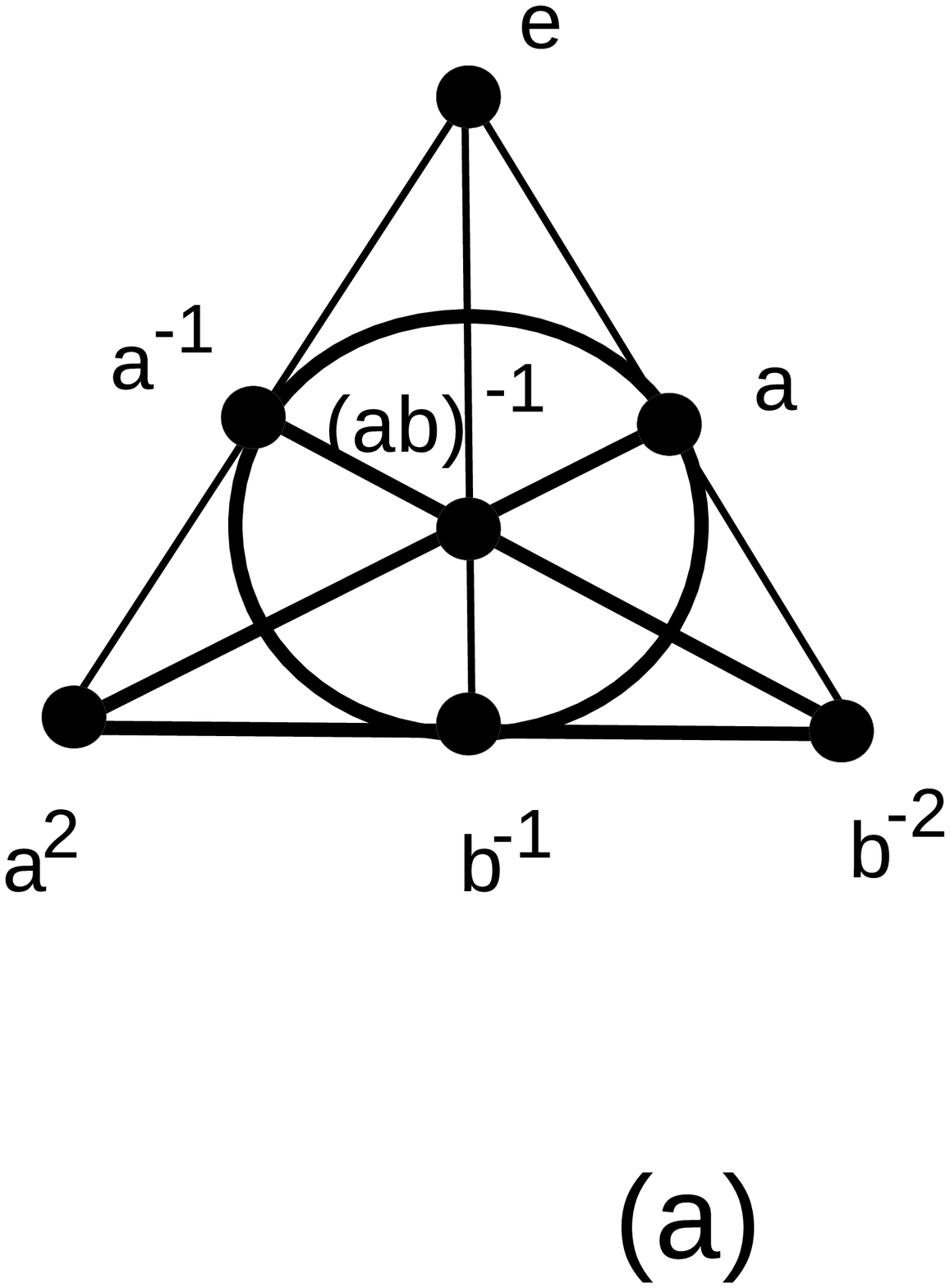}
\includegraphics[width=5cm]{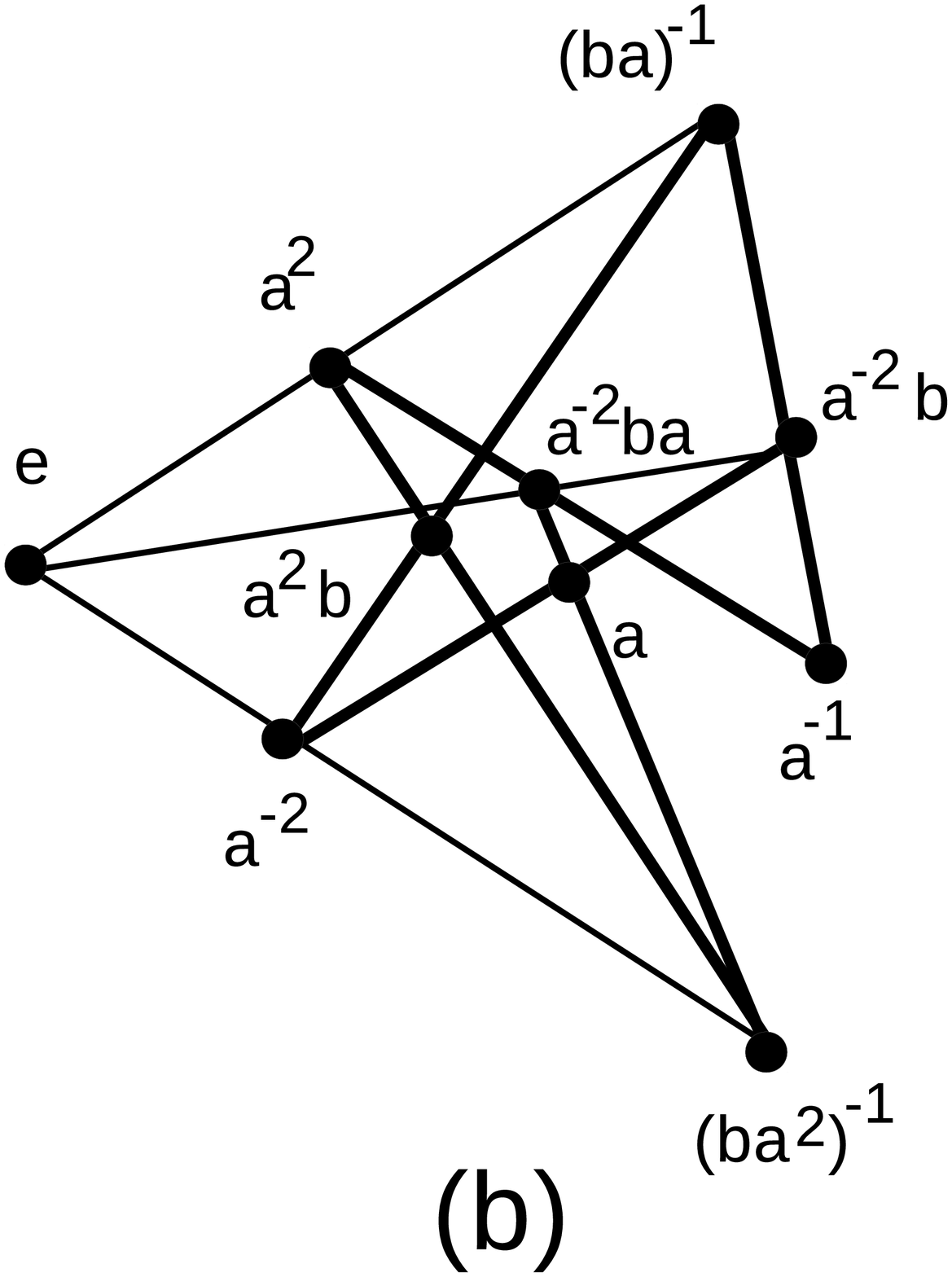}
\caption{Contextual geometries associated to \lq (i) true and (ii) false' for the MICs of the Figure-of-Eight knot $K4a1$ listed in Table \ref{tab1}: (a) the Fano plane related to the manifold otet$14_{00002}$ at index $7$, (b) the configuration $[10_3]$ at index $10$. The bold lines are for cosets that are not all mutually commuting. Each line corresponds to pair of cosets with the same stabilizer subgroup isomorphic to $\mathbb{Z}_2^2$.}
\label{fig3}
\end{figure}

There are three exceptions where (i) is true and a geometry does exist (when (ii) fails to be satisfied). This first occurs in dimension $4$ with a $2QB$ MIC arising from the $3$-manifold otet$08_{00002}$, in this case geometric contextuality occurs in the cosets as in Fig. \ref{fig2} of the previous subsection. Then, it occurs in dimension $7$ (corresponding to $3$-manifolds otet$14_{00002}$ and otet$14_{00035}$) when the geometry of cosets is that of the Fano plane shown in Fig. \ref{fig3}a. Finally it occurs in dimension $10$ when the geometry of cosets is that of a $[10_3]$ configuration shown in Fig. \ref{fig3}b \cite[p 74]{Grunbaum}. 

In addition to the latter cases, false detection of a MIC may occur (this is denoted \lq fd') in dimension $8$ as shown in Table \ref{tab2}.

\small

\begin{table}[ht]
\begin{center}
\begin{tabular}{|c|l|c|c||r|l|r|}
\hline 
\hline
d & ty& $M^3$ (or $P$)& cp & (i) & pp & geometry\\
\hline
2 &cyc&  otet$04_{00002}$, $m206$   &1 &no & &\\
\hline
3 & cyc& otet$06_{00003}$,  $s961$  &1  &no& &\\
\hline
4 &irr&  otet$08_{00002}$, $L10n46$, $t_{12840}$ &2  & yes & 2 & 2QB MIC, $\bm{\Delta}$\\
  &cyc&   otet$08_{00007}$, $t12839$ &1  &no & & \\
\hline
5 & cyc&  otet$10_{00019}$ &1 & no   &   &\\
  & irr& otet$10_{00006}$, $L8a20$ &3 & yes & &$\Delta$\\
  & irr $(\times 2)$&  otet$10_{00026}$ &2 & yes & 1 & $5$-dit MIC\\
\hline 
6  &cyc& otet$12_{00013}$ &1 &no &  & \\
   &irr&   otet$12_{00039}$ &1 & no &  & \\
   &irr $(\times 2)$&   otet$12_{00038}$ &1 & yes & 10 &  $6$-dit MIC\\	
   &irr & otet$12_{00041}$ &2 & no&  & \\   
   &irr $(\times 2)$& otet$12_{00017}$  &2& no &  & \\
   &irr $(\times 4)$&  otet$12_{00000}$ &2 &yes &2  & $6$-dit MIC\\
\hline
7  &cyc& otet$14_{00019}$ & 1 & no  &  \\
   &irr $(\times 4)$& otet$14_{00002}$, $L14n55217$ &3 & yes & 2  & $7$-dit MIC,$\bm {\Delta:  \mbox{Fano}}$ \\
   &irr $(\times 4)$& otet$14_{00035}$ &1 & yes & 2  &  $7$-dit MIC, $\bm{\Delta: \mbox{Fano}}$\\
	\hline
8  &cyc& otet$16_{00026}$ & 1 & no  &  & \\
	 &irr $(\times 2)$& otet$16_{00035}$ & 1 & no & &  \\
	 &irr $(\times 2)$& otet$16_{00079}$ & 2 & yes & &  fd\\
	 &irr $(\times 2)$& otet$16_{00016}$ & 2 & yes & & fd \\
	 &irr & otet$16_{00092}$ & 2 & no & &  \\
	 &irr & otet$16_{00091}$ & 2 & yes & & $16$-cell  \\
	 &irr & otet$16_{00013},L14n17678$ & 2 & no & &  \\
		\hline
    \hline
	9  && $(36,9)\cong 3^2 \rtimes 4$ &  & yes  & 2 & 2QT MIC\\	
	  &$(\times 2)$& $(504,156)=PSL(2,8)$ &  & yes  & 3 & 2QT MIC\\
		&$(\times 2)$& $(216,153)\cong 3^2 \rtimes (24,3)$ &  & yes  & 2 & 2QT MIC\\
		\hline
		10  &$(\times 6)$& $(160,234) \cong 2^4 \rtimes 10$ &  & yes  & 5 & $10$-dit MIC\\
			  &$(\times 2)$& $(120,34)=S_5$ &  & yes  & 4 & $10$-dit MIC, $\bm{\Delta: [10_3]}$\\
					&$(\times 2)$& $(120,34)=S_5$ &  & no  & 7 & $10$-dit MIC, $5$-ortho\\	
				&& $(360,118)=A_6$ &  & yes  & 5 & $10$-dit MIC\\	
		\hline
    \hline
\end{tabular}
\caption{\footnotesize{Table of $3$-manifolds $M^3$ found from subgroups of finite index $d$ of the fundamental group $\pi_1(S^3\setminus K4a1)$ (alias the $d$-fold coverings over the Figure-of-Eight knot $3$-manifold). The covering type \lq ty' in column 2, the manifold identification \lq $M^3$' in column 3 and the number of cusps \lq cp' in column 4 are from SnapPy \cite{SnapPy}. For $d=9$ and $10$, SnapPy does not provide results so that we only identify the permutation group $P=$SmallGroup$(o,k)$ (abbreviated as $(o,k)$), where $o$ is the order and $k$ is the $k$-th group of order $o$ in the standard notation (that is used in Magma). If it exists, the MIC is \lq $pp$'-valued. If (i) is true, unless otherwise specified, the graph of cosets leading to a MIC is that of the $d$-simplex [and/or the condition (ii) is true: no geometry]. The symbol $\Delta$ means that (ii) fails to be satisfied. When there exists a MIC with (i) true and (ii) false, the geometry is shown in bold characters. The symbol \lq fd' means a false detection of a MIC when (i) and (ii) are satisfied simultaneously while a MIC does not exist. The abbreviations \lq Fano', \lq $d$-ortho' and \lq $[10_3]$' are  for the Fano plane, the $d$-orthoplex and the corresponding geometric configuration.}} 
\label{tab2}
\end{center}
\end{table}
%

\newpage

\section{MIC states pertaining to the Trefoil knot and its $0$-surgery}
\label{TrefoilKnot}

We now investigate the relation of MIC states to the group geometrical axioms (i)-(ii) (or their negation) in the context of the trefoil knot $3_1$ (in Sec. \ref{trefoilb}) and its $0$-sugery (in Sec. \ref{trefoila}).
The fundamental group of the complement of $3_1$ in the $3$-sphere $G=\pi_1(S^3 \setminus 3_1)$, and its connection to MICs, is studied in \cite[Table 1]{MPGQR1} and below.

\subsection{Group geometrical axioms applied to the fundamental group $\bm{\pi_1(\tilde{E}_8)}$}
\label{trefoila}

The manifold $\tilde{E}_8$ is defined by $0$-surgery on the trefoil knot $3_1$ and is of special interest as shown in \cite[Sec 3]{MPGQR3} and references therein. The number of subgroups of index $d$ of the fundamental group $\pi_1(Y)$ is as follows

\small
\begin{table}[h]
\begin{center}
\begin{tabular}{|c|l|c|c||r|l|r|}
\hline 
\hline
d & P & (i) & pp & geometry\\
\hline
3  & $6$ &yes &1 & Hesse SIC, $\bm{\Delta}$\\
4  & $A_4$&yes &2 & 2QB MIC, $\bm{\Delta}$\\
6  & $A_4$ &no &2 & 6-dit MIC, $K(2,2,2)$\\
7  & $(42,1)\cong 7\rtimes(6,2)$ &yes &2 & 7-dit MIC\\
9  & $(54,5)\cong 3^2 \rtimes (18,3)$, $(\times 2)$ &yes & & $K(3,3,3)$\\
12  & $(72,44) \cong 2^2 \rtimes (18,3)$  &yes & & $\overline{\mathcal{L}(K(3,4))}$\\
13  & $(78,1) \cong  13 \rtimes (6,2)$, $(\times 2)$ &yes &4 & 13-dit MIC\\
15  & $(150,6) \cong 5^2 \rtimes (6,2)$, $(\times 2)$ &no &6 & 15-dit MIC, $K(5,5,5)$\\
16  & $(96,72) \cong 2^3 \rtimes A_4$ &yes & & $K(4,4,4,4)$\\
19  & $(114,1) \cong 19 \rtimes (6,2)$ &yes &3 & 19-dit MIC\\
21  & $(126,9) \cong 7 \rtimes (18,3)$, $(\times 2)$ & yes&5 & \tiny{21-dit MIC, $\Delta$: {\bf K(3,3,3,3,3,3,3)}}\\
\hline
\hline
\end{tabular}
\caption{Table of subgroups of the fundamental group $\pi_1[S^3\setminus 3_1(0,1)]$ [with $3_1(0,1)$ the 0-surgery over the trefoil knot] when the condition (i) is satisfied or when a MIC is missed. See the captions of Table \ref{tab1} and \ref{tab2} for the meaning of abbreviations. } 
\label{tab3}
\end{center}
\end{table}

$\eta_d[\tilde{E}_8]=[1,1,{\bf 2},{\bf 2},1,~{\bf 5},{\bf 3},2,4,1,~1,12,{\bf 3},3,{\bf 4},~3, 1, 17, {\bf 3}, 2,~{\bf 8}, 1, 1, 27, 2,~\ldots]$

\noindent
where a bold number means that a MIC exists at the corresponding index.

Such cases are summarized in Table \ref{tab3}. As expected this occurs when the axioms (i) and (ii) are both true, or are both false. The latter case occurs at index $6$ with geometry of the octahedron [and graph $K(2,2,2)$] and at index $15$ with a geometry of graph $K(5,5,5)$.

Exceptions to the rules are when a MIC exists with (i) true but not (ii). This occurs in dimension $3$ (for the Hesse SIC) since the free group has a single generator (a trivial case), at index $4$ (for the $2QB$ MIC) with a contextual geometry as in Fig. \ref{fig2} and at index $21$ with a contextual geometry (not shown) of graph $K(3,3,3,3,3,3,3)$.

\subsection{Group geometrical axioms applied to the fundamental group $\bm{\pi_1(S^3 \setminus 3_1)}$}
\label{trefoilb}

The characteristics of submanifolds obtained from the subgroups of index $d$ of the fundamental group $\pi_1(S^3 \setminus 3_1)$ for the trefoil knot complement are given in Table \ref{tab4} by using SnapPy \cite{SnapPy} and Sage \cite{Sage} for identifying the corresponding subgroup of the modular group $\Gamma$ \cite{PlanatModular} (this improves \cite[Table 1]{MPGQR1}).

As for the above sections, when axioms (i) and (ii) are simultaneously satisfied (or both are not satisfied), a MIC is created. Otherwise no MIC exist in the corresponding dimension, as one should expect.

\begin{figure}[h]
\includegraphics[width=5cm]{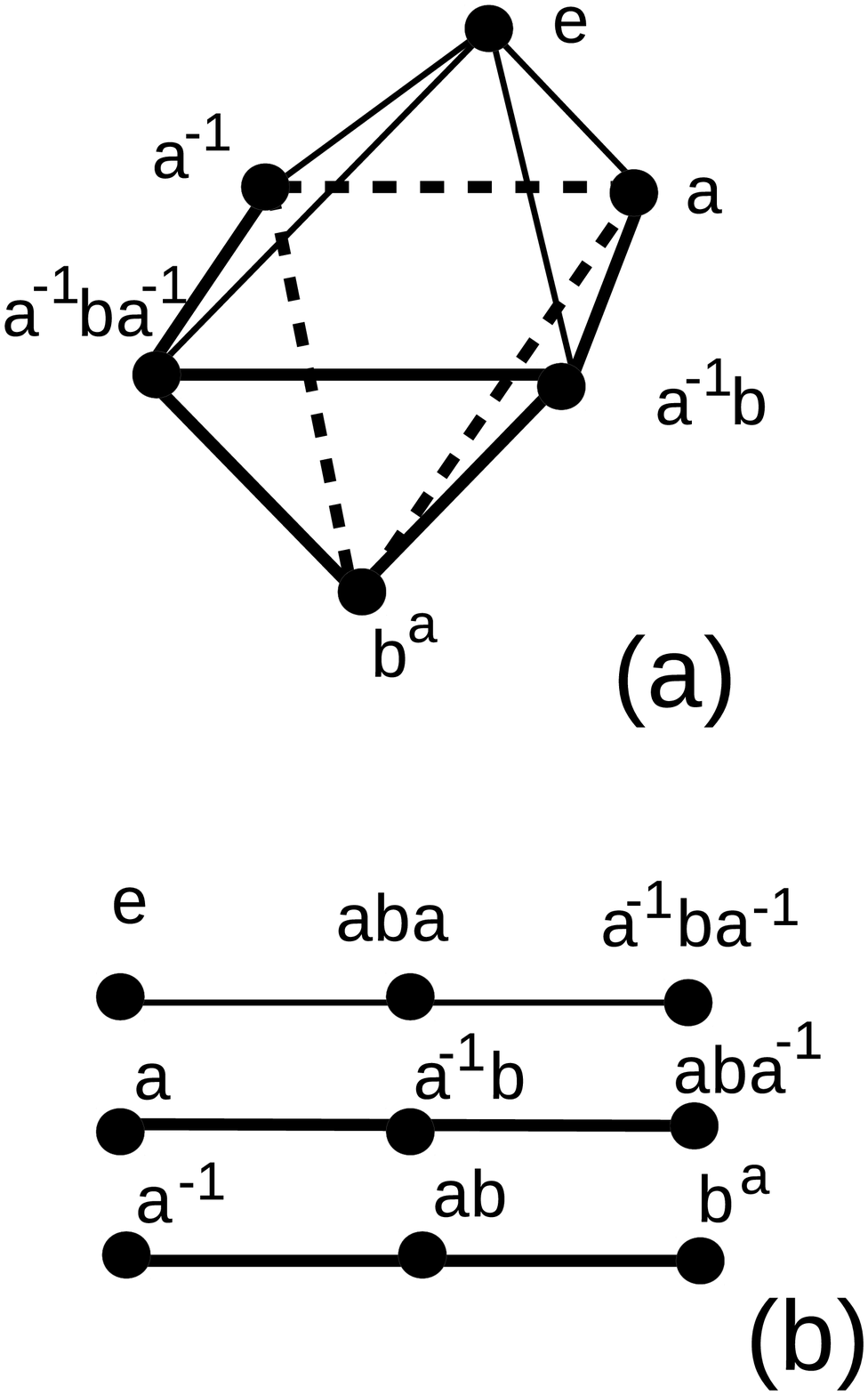}
\includegraphics[width=5cm]{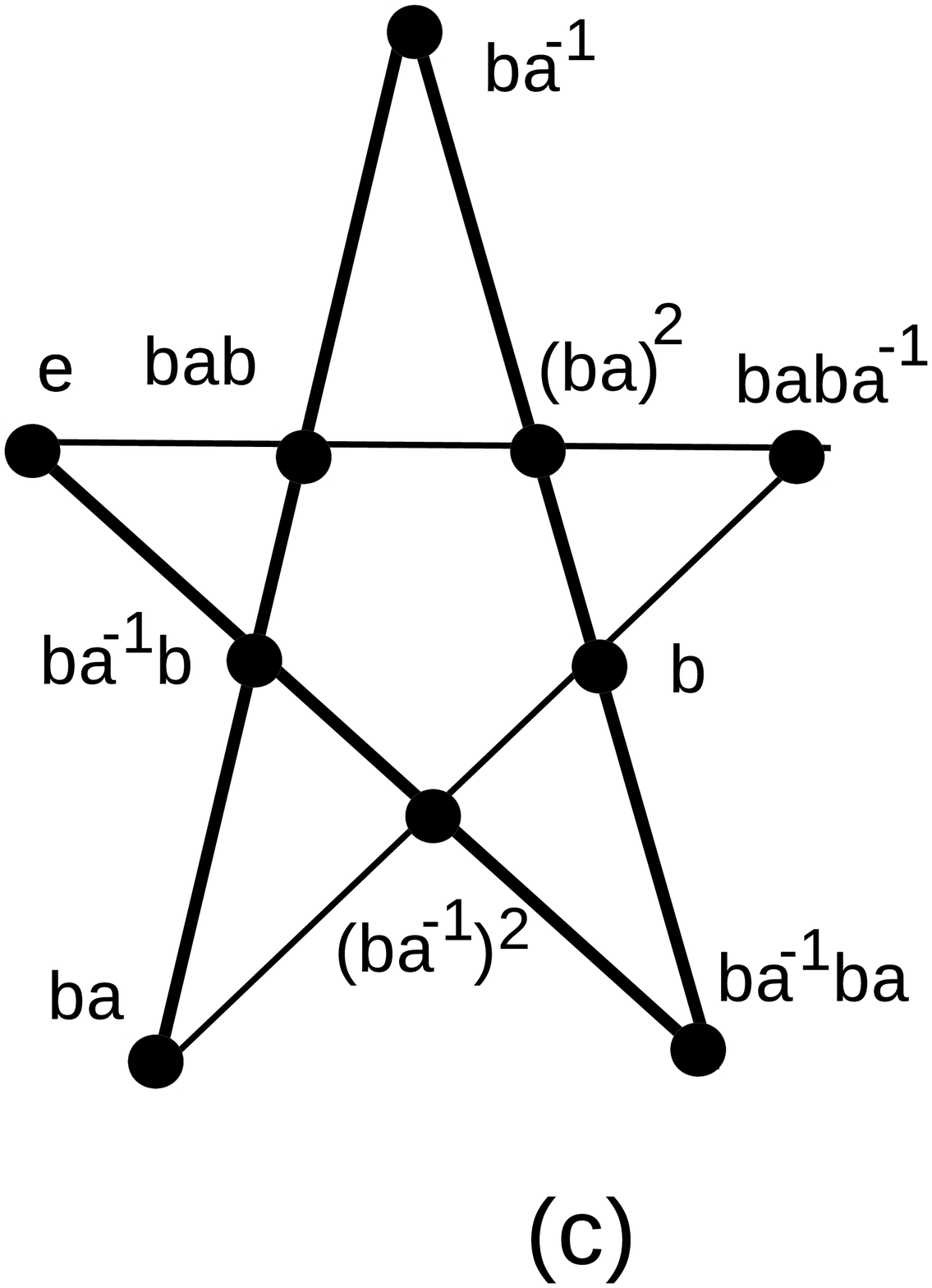}
\caption{Contextual geometries associated to \lq (i) true and (ii) false' for the MICs of the trefoil knot $3_1$ listed in Table \ref{tab2}: (a) the octahedron, related to the subgroup $\Gamma_0(4)$ of $\Gamma$ at index $6$ , (b) three disjoint lines $K_3^3$ at index $9$, (c) The Mermin's pentagram at index $10$. The bold lines are for cosets that are not all mutually commuting.}
\label{fig4}
\end{figure}

There are a few exceptions where (i) is true and a geometry does exist (when (ii) fails to be satisfied). This first occurs in dimension $3$ for the Hesse SIC where the free subgroup is trivial with a single generator. The next exceptions are for the $6$-dit MIC related to the permutation group $S_4$ with the contextual geometry of the octahedron shown in Fig. \ref{fig4}a, in dimension $9$ for the $2QT$ MIC related to the permutation group $3^3 \rtimes S_4$ with a contextual geometry consisting of three disjoint lines, and in dimension $10$ for a $10$-dit MIC related to the permutation group $A_5$ and the contextual geometry of the so-called Mermin pentagram. The latter geometry is known to allow a $3QB$ proof of the Kochen-Specker theorem \cite{Planat2012}.

\small
\begin{table}[h]
\begin{center}
\begin{tabular}{|c|l|c||c|r|l|r|l|}
\hline 
\hline
d & ty & cp & P & (i) & pp & type in $\Gamma$ & geometry\\
\hline
2 & cyc & 1 & (2,1) $\equiv$ 2 & no &&& \\
\hline
3 & cyc & 1 & (3,1) $\equiv$ 3 & no & & & \\
  & irr & 2 & (6,1) $\equiv$ 6 & yes & 1& $\Gamma_0(2)$& \scriptsize{Hesse SIC, $\bm{\Delta}$, L7n1} \\
	\hline
	4 & cyc & 1 & (4,1) $\equiv$ 4 & no & & & \\
	 & irr & 2 & $(12,3)=A_4$ & yes & 2 & $\Gamma_0(3)$ & \scriptsize{2QB MIC, $\bm{\Delta}$, L6a3} \\
	& irr & 1 & $(24,12)=S_4$ & yes & 2 & $4A^0$ & 2QB MIC \\
	\hline
5 & cyc & 2 & (5,1) $\equiv$ 5 & no & & & \\	
 & irr & 3 & $(60,5)=A_5$ & yes & 1& $5A^0$& $5$-dit MIC\\
\hline
6 & reg & 3 & (6,1) $\equiv$ 6 & no & 2 & $\Gamma(2)$ & \scriptsize{$6$-dit MIC, $6_3^3$ \cite{MPGQR3}} \\	
 & cyc & 3 & (6,2) = 3$\times 2$ & no &  & $\Gamma'$ &  \\
 & irr & 2 & $A_4$ & no & 2 & $3C^0$ & \scriptsize{$6$-dit MIC, K(2,2,2)} \\
& irr & 1 & $(24,13)= 3 \rtimes 8$ & no &  & $6B^0$ &  \\
& irr & 1 & $(18,3)\cong 3^2 \rtimes 2$ & no &  & $6A^0$ &  \\
& irr & 3 & $S_4$ & yes & 2 & $\Gamma_0(4)$ &\tiny{ $6$-dit MIC, $\bm{\Delta: \mbox{octa}}$}  \\
& irr & 2 & $A_5$ & yes & 2 & $\Gamma_0(5)$ & $6$-dit MIC\\
& irr & 2 & $S_4$ & yes & 2 & $4C^0$ & \tiny{$6$-dit MIC, $\bm{\Delta: \mbox{octa}}$}\\
\hline
7 & cyc & 1 & (7,1)$\equiv 7$ & no & & & \\
 & \scriptsize{irr ($\times 2$) }&2 & \scriptsize{$(42,1)\cong 7 \rtimes (6,2)$} & yes & 2& \tiny{$NC(0,6,1,1,[1^16^1])$} & $7$-dit MIC\\
& \scriptsize{irr $(\times 2)$ }& 1 & \scriptsize{$(168,42)=PSL(2,7)$} & yes & 2 & $7A^0$ & $7$-dit MIC\\
& \scriptsize{irr ($\times 2$)} &2 & $S_7$ (order 5040) & yes && \tiny{$NC(0,10,1,1,[2^15^1])$} & ? $7$-dit MIC\\
\hline
8 & cyc & 1 & $(8,1)\equiv 8$ & no & & & \\
 & irr & 2 & (24,13) & no & &$6C^0$ & \\
& irr & 2 & $S_4$ & no & &$4D^0$ & \\
&\scriptsize{irr ($\times 2$)}& 2 & $(24,3)\cong 2. A_4$ & yes & & &  $16$-cell \\
& irr  & 2 & $PSL(2,7)$ & yes & &$\Gamma_0(7)$ &  fd \\
& \scriptsize{irr ($\times 2$) } & 1 & $SL(2,7)$ & yes & &\tiny{$NC(0,8,2,2,[8^1])$} &  fd \\
& \scriptsize{irr ($\times 2$) }& 2 & \scriptsize{$(48,29)\cong 2.(24,3)$} & yes & &$8A^0$ & $16$-cell  \\
\hline 
\hline
9 &  &  & $(9,1)\equiv 9$ & no & & & \\
 &  & 2 & (18,3) & no & 7 & $6D^0$ & \scriptsize{$9$-dit MIC, K(3,3,3)} \\
&  & 2 & \scriptsize{$(54,5) \cong 3^2 \rtimes (18,3)$} & no & 7 & \tiny{$NC(0,6,3,0,[3^16^1])$} &\scriptsize{ $9$-dit MIC, K(3,3,3)} \\
&  & 1 & \scriptsize{$(324,160)\cong 3^3 \rtimes A_4$} & no & &$9A^0$ &  K(3,3,3), $K_3^3$ \\
&  & 3 & (54,5) & yes & 7 & \tiny{$NC(0,6,1,0,[1^12^16^1])$}  & $9$-dit MIC \\
 & ($\times 3$) &  &\scriptsize{ $(162,10)\cong 3^2 \rtimes 6$ }& yes & & & K(3,3,3) \\
& ($\times 2$) & 1 & \scriptsize{$(504,156)=PSL(2,8)$} & yes & 3 &\tiny{$NC(1,9,1,0,[9^1])$}&2QT MIC \\
& ($\times 2$) & 2 & \tiny{$(432,734)\cong 3^2 \rtimes (48,29)$} & yes & 2 &\tiny{$NC(0,8,3,0,[8^11^1])$}&2QT MIC \\
&  & 3 & \scriptsize{$(648,703)\cong 3^3 \rtimes S_4$} & yes & 2 &\tiny{$NC(0,12,1,0,[2^13^14^1])$}&\tiny{2QT MIC, $\bm{\Delta: K_3^3}$} \\
\hline
10&  & 1 & \scriptsize{$(120,35) \cong 2 \rtimes A_5$ }& no & &$10A^0$& $5$-ortho \\
& &   2& $A_5$ & yes & 5 &$5C^0$& $10$-dit MIC,\scriptsize{ $\Delta$: {\bf MP} }\\
& ($\times 2$) &  1& \scriptsize{$(720,764) \cong A_6 \rtimes 2$ }& yes & 9 &\tiny{$NC(0,10,0,4,[10^1])$}& $10$-dit MIC\\
\hline 
\hline
\end{tabular}
\caption{ \footnotesize{Subgroups of index $d$ of the fundamental group $\pi_1(S^3\setminus 3_1)$ (alias the $d$-fold coverings over the trefoil knot $3$-manifold. The meaning of symbols is as in Table 2. When the subgroup in question is a subgroup of the modular group $\Gamma$, it is identified as a congruence subgroup or by its signature $NC(g,N,\nu_2,\nu_3,[c_i^{W_i}])$ (see \cite{PlanatModular} for the meaning of entries). The permutation group $P=$SmallGroup$(o,k)$ is abbreviated as $(o,k)$. As in Table 1, if (i) is true, unless otherwise specified, the graph of cosets leading to a MIC is the $d$-simplex [and/or the condition (ii) is true].  Exceptions (with geometry identified in bold characters) are for a MIC with (i) true and (ii) false.  For index $9$ and $10$, some subgroups of large order could not be checked as leading to a MIC or not, they are not shown in the table. The abbreviation \lq octa' is for the octahedron,  \lq MP' is for the Mermin pentagram and $K_3^3$ means three disjoint triangles.}} 
\label{tab4}
\end{center}
\end{table}
\normalsize

\section{Conclusion}

Previous work about the relationship between quantum commutation and coset-generated finite geometries has been expanded here by establishing a connection between coset-generated magic states and coset-generated finite geometries. The magic states under question are those leading to MICs (with minimal complete quantum information in them). We found that, given an appropriate free group $G$, two axioms (i): the normal closure $N$ of the subgroup of $G$ generating the MIC is $G$ itself and (ii): no coset-geometry should exist, or the negation of both axioms (i) and (ii), are almost enough to classify the MIC states. The few exceptions rely on configurations that admit geometric contextuality. We restricted the application of the theory to the fundamental group of the $3$-manifolds defined from the Figure-of-Eight knot (an hyperbolic manifold) and from the trefoil knot, and to $4$-manifolds $Y$ and $\tilde{E}_8$ obtained by $0$-surgery on them. It is of importance to improve of knowledge of the magic states due to their application to quantum computing and we intend to pursue this research in future work.

\end{document}